\input amstex
\input optlab

\def\half{{\tfrac12}}

\def\v{\vert}
\def\eps{\epsilon}
\def\varepsilon{\epsilon}
\hfuzz=3mm
\def\wt{\widetilde}

\input epsf
\documentstyle{grgppt}
\redefine\qed{\hfill$\square$}
\pagewidth{125mm}
\pageheight{185mm}
\parindent=8mm
%\TagsOnRight
\magnification=1200
\def\section#1{\bigskip\goodbreak\line{{\bf\ignorespaces #1}\hfill}
  \nobreak\medskip\nobreak\noindent\ignorespaces}

\def\subsection#1{\goodbreak\medskip
  \flushpar{\smc #1}\nobreak\smallskip \nobreak\flushpar\ignorespaces}
\def\capt#1#2{\baselineskip=10pt\botcaption{\baselineskip=10pt\eightpoint{\it Figure #1.}\q
  #2}\endcaption\par}
\def\mletter#1#2#3{\hskip#2cm\lower#3cm\rlap{$#1$}\hskip-#2cm}
\def\lastletter#1#2#3{\hskip#2cm\lower#3cm\rlap{$#1$}\hskip-#2cm\vskip-#3cm}
\def\figure#1\par{\parindent=0pt
  \vbox{\baselineskip=0pt \lineskip=0pt
  \line{\hfil}
  #1}}
\def\Bigmid{\,\Big\vert\,}
\def\BG{[\ci{BG1}]}
\def\i{\item}
\def\bzero{0}
\def\bx{\bold x}
\def\bp{\bold p}
\def\bv{\bold v}
\def\bu{\bold u}
\def\boldy{\bold y}
\def\q{\quad}

\def\S{{\text{\rm s}}}
\def\T{{\text{\rm t}}}

\def\la{\langle}
\def\ra{\rangle}
\def\b{\beta}

\def\a{\alpha}
\def\t{\theta}

\def\g{\gamma}
\def\sF{\Cal F}

\def\oo{\infty}

\def\L{{\Bbb L}}
\def\Z{{\Bbb Z}}
\def\E{{\Bbb E}}
\def\R{{\Bbb R}}
\def\P{{\Bbb P}}
\def\ldc{\vec\L^d}
\def\edc{\vec\E^d}
\def\ldca{\vec\L^d_{\text{alt}}}
\def\edca{\vec\E^d_{\text{alt}}}
\def\ltwoc{\vec\L^2}
\def\lthreeca{\vec\L^3_{\text{alt}}}

\def\ltwoca{\vec{\L}^2_{\text{alt}}}

\def\pc{{p_{\text{\rm c}}}}
\def\wt#1{\widetilde #1}
\def\ci#1{\citeto{#1}}

\def\rc{r_{\text{c}}} 
\def\fc{\to\hskip-3pt{}_{\text{fc}}\hskip3pt}

\catcode`\@=11
 \def\logo@{}
\catcode`\@=13
\def\today{\number\day
\space\ifcase\month\or
  January\or February\or March\or April\or May\or June\or
  July\or August\or September\or October\or November\or December\fi
  \space\number\year}
\topmatter
\title Directed percolation and random walk
\endtitle
\author Geoffrey Grimmett and Philipp Hiemer
\endauthor
\address
Statistical Laboratory, University of Cambridge, Wilberforce Road,
Cambridge CB3 0WB, United Kingdom
\endaddress
\email g.r.grimmett{\@}statslab.cam.ac.uk, p.hiemer{\@}statslab.cam.ac.uk\endemail
\http http://www.statslab.cam.ac.uk/$\sim$grg/ \endhttp

\abstract
Techniques of `dynamic renormalization', developed earlier
for undirected percolation and the contact model, 
are adapted to the setting of directed percolation, thereby obtaining
solutions of several problems for directed percolation
on $\Z^d$ where $d \ge 2$.
The first new result is a type of uniqueness theorem:
for every pair $x$ and $y$ of vertices which lie in infinite
open paths, there exists almost surely
a third vertex $z$ which is joined to infinity and which
is attainable from $x$ and $y$ along directed open paths.
Secondly, it is proved that a random walk on an infinite
directed cluster is  transient, almost surely, when $d \ge 3$.
And finally, the block arguments of the paper may be adapted
to systems with infinite range, subject to
certain conditions on the edge probabilities.
\endabstract
\endtopmatter\footnote""{This version was prepared on 6 March 2001.}

\document
\section{1. Introduction}
`Block arguments' constitute a fundamental technique
for studying disordered spatial processes. For many years,
physicists have appealed to the theory of renormalization,
although difficulties have emerged in making such arguments 
mathematically rigorous. Block arguments have been
widely used since [\ci{Kes81}], at least,
and have led to proofs of several theorems of
substance (see [\ci{BGN1}, \ci{BGN2}, \ci{Cerf0}, 
\ci{GrM}, \ci{Pis95}]). For a recent
account of dynamic and static renormalization 
in the context of undirected percolation,
the reader is referred to Chapter 7 of [\ci{G99}].

The adaptation of such methodology to directed
models, such as directed percolation, is not totally straightforward.
New difficulties arise
through the presence of orientations on the edges of the underlying
lattice, and new ideas are needed to overcome these problems.
Various specific problems have arisen in the work of other authors on
directed percolation. These problems
may be dealt with by suitable block arguments,
and the target of this paper is to show how this may be done.

We introduce  directed percolation in Section 2, and we state
our main results in Section 3. We proceed
in Section 4 to describe the block construction
which enables a comparison between a
supercritical process and another directed
percolation process with density close to $1$.  In the manner
already explored in [\ci{BGN2}, \ci{BG1}], 
this implies the absence of infinite open paths
in the critical case, as well as the continuity of the percolation probability.

As consequences of the above
comparison theorem, we shall present positive answers 
to the following two questions posed
by others. Let $d\ge 2$,  let $\pc$ be the critical probability
of directed percolation on $\Z^d$, and suppose $p>\pc$. 
We say that a vertex $x$ 
is {\it connected to infinity\/}
if there exists an infinite open path which is
directed away from its endvertex $x$.
It is proved that, almost surely,
for all vertices $x$, $y$ which are connected to infinity, there exists a third
vertex $z$, also connected to infinity, such that $x$ and $y$ are connected
by directed open paths to $z$. This answers a question of 
Itai Benjamini, and is a
cousin of the `unique infinite cluster' theorem of
[\ci{AKNa},  \ci{BK}, \ci{GanGR}].

Our second application concerns the classification of a random
walk on the set of vertices attainable from the origin along directed
open paths.  It was proved in [\ci{GKZ}] that random
walk on the infinite open cluster of undirected percolation is
almost surely transient in three or more dimensions. A similar
result was proved in [\ci{BenPP97}]
for directed percolation in three  dimensions, whenever the edge density is
sufficiently large. The methods used in the latter paper are quite different
 from
those used in [\ci{GKZ}], and have other applications also.
We present in Theorem 3 a positive answer to a question posed in 
[\ci{BenPP97}], namely whether the transience result for
directed percolation may be extended to all values of $p$ satisfying $p>\pc$.

There is only little difficulty in extending results
about `nearest neighbour' directed percolation
to systems with {\it finite\/} range. There has been 
recent interest [\ci{vdHS00}, \ci{Sl00}]
in systems with {\it infinite\/} range, and particularly in whether
or not they may possess an infinite open path {\it at\/}
the critical point. Such systems are explored in Section 7, 
where it is explained 
how the block construction of Section 4 may be adapted
to infinite-range systems satisfying some weak conditions of regularity.

Directed percolation is closely related to the contact model, for
which block arguments have been used to prove results related to
some of those described above 
(see [\ci{BG1}, \ci{Dur91}, \ci{Lig}, \ci{Lig99}]).
Although the arguments of [\ci{BG1}] are in part useful for
the present work, the discreteness of the underlying lattice  
leads to some special problems for directed percolation. Just as in the case
of the contact model, the comparison theorem of the current paper may be used
to establish further results such as a shape theorem, a complete 
convergence theorem, and the continuity of the critical points of slabs
in the limit of large slabs. 
We do not present the details of the necessary
proofs; an interested reader may refer to the earlier papers
cited in [\ci{BG1}], where closely related material is studied.

\section{2. Notation}
Let $\Z^d$ denote the set of all $d$-vectors $x=(x_1,x_2,\dots,x_d)$
of integers. 
For $x,y\in \Z^d$, we define
$$
|x-y|=\sum_{i=1}^d|x_i-y_i|.
$$
We refer to vectors in $\Z^d$ as {\it vertices\/}, and
we turn $\Z^d$ into a graph by adding an (undirected) edge $\la x,y\ra$ between
every pair $x$, $y$ of vertices such that $|x-y|=1$. The resulting
graph is denoted $\L^d=(\Z^d,\E^d)$. The {\it origin\/} of this graph 
is the vertex $\bzero=(0,0,\dots,0)$. We write $x\le y$ if $x_i\le y_i$
for $1\le i\le d$.
We may use the lattice $\L^d$ to generate a multiplicity of directed graphs
of which two feature in this paper.

\smallskip
\flushpar{\bf Conventional model.} The edge $\la x,y\ra$ 
with $x\le y$ is assigned
an arrow from $x$ to $y$. We write $\ldc=(\Z^d,\edc)$ 
for the ensuing directed graph.

\smallskip
\flushpar{\bf Alternative model.} Each vertex $x=(x_1,x_2,\dots,x_d)$ may be
 expressed
as $x=(\bx,t)$ where $\bx=(x_1,x_2,\dots,x_{d-1})$ and $t=x_d$.
Consider the directed graph with vertex set $\Z^d$ and with a directed
edge joining two vertices $x=(\bx,t)$ and $y=(\boldy,u)$ whenever
$\sum_{i=1}^{d-1}\v y_i-x_i\v \le 1$ and $u=t+1$. We write
$\ldca=(\Z^d,\edca)$ for the ensuing directed graph, and we note 
that every vertex has out-degree $1+2(d-1)$.
\smallskip

\topinsert
\figure
\centerline{\epsfxsize=7cm\epsffile{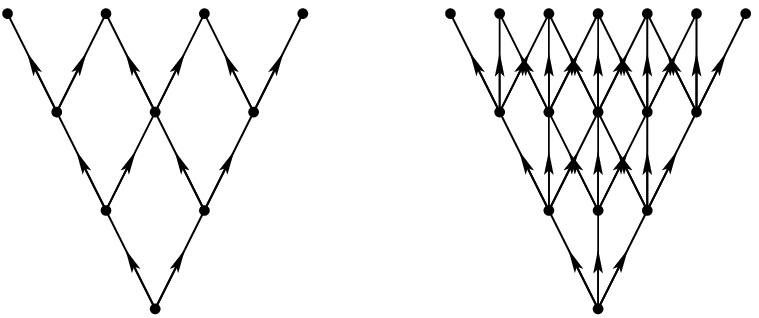}}

\capt{1}{The graphs $\ltwoc$ and $\ltwoca$.}
\endinsert

The graphs $\ltwoc$ and $\ltwoca$ are sketched in Figure 1.
 Until recently, the conventional
model has been considered the natural habitat of directed percolation
in $d$ dimensions, whereas recent results of van der Hofstad and
Slade [\ci{vdHS00}, \ci{Sl00}]
concerning the scaling limit of critical directed 
percolation in high dimensions have indicated
the relevance of the alternative model. It is reasonable to think that
results available for either of 
these models may be derived for the other also, 
but some technical difficulties may arise in justifying 
this statement in concrete examples, owing to
the fact that the automorphism groups of $\ldc$ and $\ldca$
are different. We shall in this paper 
concentrate on the alternative model, since we 
avoid thus certain minor complications involving periodicity. Our arguments
are easily adapted to the 
conventional
model. In either case, we write $[x,y\ra$ for
an edge which is directed from $x$ to $y$. 
{\it We assume henceforth that we are studying
the directed graph $\ldca$ where $d\ge 2$.}

In Section 7, we shall consider a generalization
of the alternative model in which random directed edges of long range are 
added to $\Z^d$, rather then merely between nearest neighbours of $\ldca$.
It turns out that, subject to certain natural assumptions on the
parameters of such a process, the techniques developed below and elsewhere
may be adapted successfully to such a model. Such results
have potential applications to the work reported in [\ci{vdHS00}, \ci{Sl00}], where
it is proved that the scaling limit of critical directed percolation
is, for high dimension, the process known as 
integrated super-Brownian excursion.

Given a directed or undirected graph $G=(V,E)$, 
the configuration space for percolation 
on $G$ is the set
$\Omega=\{0,1\}^{E}$. For $\omega\in\Omega$, we call an edge
$e\in E$ {\it open\/} if $\omega(e)=1$
and {\it closed\/} otherwise. With $\Omega$ we associate the
$\sigma$-field $\sF$ of subsets generated by the finite-dimensional cylinders.
For $0\le p\le 1$, we let $\P_p$ be product measure
on $(\Omega,\sF)$ with density $p$.
There is a natural partial order on $\Omega$ given as follows: for
$\omega_1,\omega_2\in\Omega$, we write $\omega_1\le\omega_2$
if $\omega_1(e)\le\omega_2(e)$ for all $e\in E$.
We shall consider primarily percolation on the graph $\ldca$,
and we suppose henceforth that $\Omega=\{0,1\}^{\edca}$.

Let $\omega \in \Omega$. An {\it open path\/} is an alternating 
sequence $x_0,e_0,x_1,e_1,x_2,\dots$
of distinct vertices $x_i$ and open edges $e_j$ such 
that $e_i=[x_i,x_{i+1}\ra$
for all $i$.  If the path is finite, it has two endvertices $x_0$,
$x_n$, and is said to connect $x_0$ to $x_n$. If  the
path is infinite, it is said
to connect $x_0$ to infinity.
A vertex $x$ is
said to be {\it connected\/} to a vertex $y$, written $x \rightarrow y$, 
if there exists an open path connecting $x$ to $y$. 
If $S\subseteq\Z^d$, we write $x \rightarrow y$ in $S$ if
there exists an open path from $x$ to $y$ using only vertices contained in $S$.
For $A,B\subseteq\Z^d $, we say that 
$A$ is {\it connected to\/} $B$ if there exist $a\in A$ and $b\in B$
such that $a\to b$; in this case, we write $A\to B$. 
We say that $A$ is {\it fully
connected to\/} $B$ if, for all
$b\in B$, there exists $a \in A$ such
that $a \rightarrow b$; in this case, we write $A\fc B$. 
Note that, if $A\fc B$ in $C$, then it is necessarily the case that
$B\subseteq C$.
If a
vertex $x$ is connected to infinity, we write $x
\rightarrow \infty$.  

For $x\in\Z^d$ and $\omega\in\Omega$, we write 
$$
C_x=C_x(\omega)= \{y \in\Z^d  : x \rightarrow y\},
$$
and we abbreviate $C_0$ to $C$. The set $C_x$ is called the
{\it open cluster at $x$}.
The
{\it percolation probability\/} is defined as the function
$$
\theta(p) = \P_p(\bzero\rightarrow\infty) = \P_p(|C|=\infty).
$$
Let
$\psi(p)=\P_p(x\rightarrow\infty\text{ for some }x\in\Z^d )$.
It is a consequence of the zero--one law that
$\psi(p)$ takes values $0$ and $1$ only, and that
$\theta(p) > 0\text{ if and only if }\psi(p)=1$;
cf.\ [\ci{G99}], Theorem (1.11). We define the {\it critical
probability\/}
$$
\pc=\sup\{p:\theta(p)=0\}.
$$

The cluster $C_x$ has been defined as the set of vertices to which $x$ is
connected. We shall at some point want to think of $C_x$ as a
directed  graph rather than
just a set of vertices, and this is achieved by adding to $C_x$
all open directed edges having both endvertices in $C_x$. The resulting
graph is denoted $\vec C_x$. The undirected graph obtained from 
$\vec C_x$ by deleting the orientations is
denoted $C_x$, and it will be clear from the context whether
$C_x$ is to be interpreted as a set of vertices or as a graph. 

\section{3. Principal results}
Our first principal result is a 
re-affirmation of a theorem of [\ci{BG1}]. The latter
paper studied the contact model rather than directed percolation,
but included some remarks on the extension of the 
results therein to directed percolation.

\proclaim{Theorem 1} Let $d \ge 2$. We have that
$\theta(\pc)=0$.  
\endproclaim

The corresponding fact for undirected percolation was proved
in [\ci{BGN1}, \ci{BGN2}, \ci{GrM}], and for the contact model in [\ci{BG1}].

One of the famous theorems of undirected percolation is the
statement that the infinite open cluster, when it exists, 
is almost surely unique (see [\ci{AKNa}, \ci{BK}, \ci{GanGR}]). 
There follows a version of 
this result for directed percolation, in answer to a question
posed in a personal communication by Itai Benjamini.
In order to state this in sufficient generality for later use, we introduce
the usual coupling of processes for different values of $p$
(see [\ci{G99}], p.\ 11). Let $\{U_e: e\in \edca\}$ be independent
random variables with the uniform distribution on $[0,1]$.
A realization of the $U_e$ is a vector $\eta\in [0,1]^{\edca}$, and we define
$$
\eta_p(e)=\cases 1&\text{if } \eta(e)<p,\\
0 &\text{otherwise}.
\endcases
$$
We call the edge
 $e$ $p${\it -open\/} if $\eta_p(e)=1$, and $p${\it -closed\/}
otherwise.
We write $\overset p \to \to$ for the relation $\to$ applied 
to the configuration $\eta_p$ (that is, for example, 
$x\overset p\to\to y$ if there exists a directed
path in $\eta_p$ from $x$ to $y$).

\proclaim{Theorem 2} Let $d \ge 2$. 
\flushpar{\rm(a)} Let $\pc<\a\le\b\le 1$. For all $x,y\in\Z^d$,
$$
\P\Bigl(\text{\rm$\exists z$ such that $x\overset\a\to\to z\overset\a\to\to\oo$
and $y\overset\b\to\to z$}\Bigmid x\overset\a\to\to\oo,\, y
\overset\b\to\to\oo\Bigr)=1.
$$
{\rm(b)} The function $\t$ is continuous on the interval $[0,1]$.
\endproclaim

Part (a) is reminiscent of results of [\ci{HP99}, \ci{Sch99}],
and part (b) of [\ci{BG1}, \ci{GrM}]. One may obtain a quantification
of part (a) which includes a lower bound on
the probability that such a $z$ exists within a given distance of $x$ and $y$,
but we do not pursue this here.

We turn now to random walk.
Let $G$ be a countably infinite connected graph with finite
vertex degrees, and let $\bzero$ be a specified vertex of $G$.
We assume for the sake of definiteness that $G$ has neither loops nor
multiple edges. Consider a random walk on the vertex set of $G$,
that is, a sequence $X_0,X_1,\dots$ of vertices such that, for each $n$,
$X_{n+1}$ is chosen uniformly at random from the neighbours
of $X_n$, each such choice being independent of all earlier choices.
Since $G$ is connected, the recurrence/transience of the random walk does
not depend on the choice of initial vertex $X_0$. We say
the $G$ is {\it transient\/} if the random walk is transient,
and we call $G$ {\it recurrent\/} otherwise.  
Initiated by the results of paper [\ci{GKZ}],
several authors have considered the question of whether or not an 
infinite open graph generated by a three-dimensional
percolation model is almost surely transient. Results
for undirected percolation include [\ci{GKZ}, \ci{HM98}, \ci{HoM}],
and the directed case has been studied in [\ci{BenPP97}] using the method of
`unpredictable paths'. Our third theorem
answers a question posed in this last paper.

\proclaim{Theorem 3} Let $d = 3$, and let $p>\pc$.
On the event $\{\v C\v=\oo\}$, the undirected
graph $C$ is almost surely transient.
\endproclaim

It is a near triviality that $C$ is recurrent 
in the corresponding statement for
two dimensions, since two-dimensional lattices are recurrent graphs
(this well-known fact is a consequence, for example,
of the results in Section 8.4 of [\ci{DoyS}] or Lemma 7.5
of [\ci{Soa}]).
It is the case that transience holds in
all dimensions $d\ge 3$; the proof of this would be similar,
and is not included here.

The remaining sections contain proofs of these theorems, followed in
Section 7 by a discussion of long- and infinite-range systems. Proofs are not
 given in their entirety, since this would be unduly long and would involve
 a considerable amount of duplication of material already published in 
[\ci{BGN1}, \ci{BGN2}, \ci{BG1}, \ci{GrM}].
Instead, we include only the extra arguments necessary for the
present setting.

\section{4. The block construction, and proof of Theorem 1}
A rigorous renormalization is the principal method introduced in  
[\ci {BGN1}, \ci {BGN2}] and developed further in [\ci{BG1}, \ci {GrM}]. The methods of
these papers may be adapted to directed percolation more or less as indicated 
explicitly in [\ci{BG1}], and we summarise this in this section.  Full details are
 omitted, since this would involve a considerable duplication of material;
 the reader is referred to \BG\ at salient points. For clarity of exposition,
 we assume throughout
that $d=3$.  The case $d=2$ is easier, proceeding by
path-intersection properties not valid in higher dimensions,
and the more general case $d \geq 3$
may be treated by extending the current notation as described in [\ci{BG1}].

Let $K$ and $L$ be positive integers, and write
$$
B(L) = [-L,L]^2 \cap \Z^2, \quad \partial B(L)=B(L)\setminus B(L-1).
$$
We refer to a box $B_{L,K}=B(L)\times [0,K]$ as a `space--time box' 
of $\Z^3$ in which $B(L)$ 
plays the role of space, and the 
final component  in $[0,K]$ plays the role of time.

The region $B_{L,K}$ has a {\it top\/} and {\it sides\/}
 given respectively as $B(L) \times \{K\}$ and $\partial B(L) \times [0,K]$.
  The top may be expressed as the union of  four squares of side length $L$
 and indexed in some arbitrary manner
with the set $\{-1, +1\}^2$.  The sides of $B_{L,K}$ are the
 union of four `facets' each of which is the union of two rectangles of
 side-length $L$ and height $K$.  We index these ensuing sub-facets 
in some arbitrary manner with the set
 $\{-1, +1\}^3$.

Let $r\geq 1$, and let $D_r=[-r,r]^2 \times \{0\}$, a `disk' centred at
 the origin.  Any translate of $D_r$  is termed an {\it $r$-disk\/}.
Let $N^{\bu}_\T(L,K)$ be the number of vertices $x$ in the
 subsquare of the top of $B_{L,K}$ indexed $\bu$ for which $D_r\to x$ in
 $B_{L,K}$. Let $N^{\bv}_\S(L,K)$ be the number of vertices $x$ in the
 sub-facet of the sides of $B_{L,K}$ indexed $\bv$ such that $D_r\to x$
 in $B_{L,K}$. The subscripts `t' and `s' stand for `top' and `sides'.

Suppose that $p$ is such that $\theta (p)>0$, and let $\varepsilon >0$. 
By a standard argument (see \BG, p. 1470), there exists an integer $r$ such that
$$
\P_p (D_r\to \infty)> 1 - \half \varepsilon^{12}\tag 4.1
$$
and we fix this value of $r$ henceforth. Cf.\ \BG,
equation (6).

Let $\a$ be the minimum of: (i) the probability that $0$ is fully connected
 to $D_r+re_3$, and (ii) the probability that $0$ is fully connected
 to $D_r+re_1 +2re_3$  by  paths using edges contained in $(D_r+re_1)
 + [0,2r]e_3$; here, $e_i$ denotes a unit vector of the lattice in the $i$th direction. 
 Let $M$ be large enough to ensure that in $M$ or more independent trials of
 an experiment having success probability $\alpha$, the probability of
 obtaining at least one success exceeds $1-\varepsilon$.  Let $N$ be large
 enough to ensure that, in any subset of $\Z^3$ having size $N$ or larger,
 there exist at least $M$ points all pairs of which are $L^\infty$-distance
 at least $3r+1$ apart.

There follows the main lemma; cf. Lemma (7) of \BG.

\proclaim{Lemma 4.1}
There exist positive integers $L,K$ such that, for every index $\bu \in \{-1,+1\}^2$
 and every index $\bv \in \{-1, +1\}^3$,
$$
\P_p \big(N^{\bu}_\T (L,K) \geq N\big) \geq  1-\varepsilon, 
\quad \P_p \big(N^{\bv}_\S
 (L,K) \geq N\big) \geq 1-\varepsilon .
$$
\endproclaim

\demo{Proof}
This follows very closely that of Lemma (7) of \BG, and we omit almost all 
details.  The only complication arises as remarked on p.\ 1473 of \BG.  Pick $R$
 sufficiently large that, in $RN$ independent trials with success probability
 $p$, the chance of at least $N$ successes exceeds $1- \tfrac{1}{8}
 \varepsilon^4$.  We now follow the argument of \BG\  with the difference that,
 instead of requiring that the number $N_\T(L,K)$ of points on the top of
 $B_{L,K}$ which are joined to $D_r$ by directed paths of $B_{L,K}$ satisfies
 $N_\T(L,K) \geq 4N$, we require instead that $N_\T(L,K) \geq 4RN$. 
 We derive the corresponding version of \BG, equation (12),
 with $1-\varepsilon^4$ replaced by $1-\tfrac{7}{8} \varepsilon^4$. 
 As in \BG, p.\ 1473,
we choose $S=S(L,M)$ such that
$$
\P_p \big(N_\T(L,S) \geq 4RN\big)< 1-\tfrac{7}{8}\varepsilon^4 \leq 
\P_p \big(N_\T(L,S-1) \geq 4RN\big).\tag 4.2
$$
By the choice of $R$, we have that
$$\align
\P_p \big(N_\T(L,S) \geq 4N\big) & \geq (1-\tfrac{7}{8}\varepsilon^4) 
\P_p \big(N_\T(L,S-1) \geq 4RN\big) \\
& \geq (1-\tfrac{7}{8}\varepsilon^4)(1-\tfrac{1}{8}\varepsilon^4)\geq 1-  
\varepsilon^4.
\endalign
$$
We now follow the argument of \BG, using the left inequality of 
equation (4.2), in
 order to obtain an inequality corresponding to equation (16) of \BG.

Note that the proof of the step corresponding to equation (15) of \BG\ is
 easier in the current setting, owing to the discreteness of
 the time variable. 
\qed\enddemo

Theorem 1 may now be proved exactly as was Theorem (1) of \BG.  The idea is
 to use the block $B_{L,K}$ of Lemma 4.1 and to iterate the construction therein in
 order to build, with large probability, a directed path within a certain tube of
 $\Z^3$.  
 As described in \BG, this enables a stochastic comparison with a certain 1-dependent
 percolation model with density which may be made close to $1$ by an
 appropriate choice of $\varepsilon$.  Since the events in Lemma 4.1
depend on the states of only a finite number of edges, their probabilities
 are continuous functions of $p$. It follows that the 
resulting block construction is 
infinite with strictly positive probability, 
for some $p'$ satisfying $p'<p$. 
 Thus, if $\theta(\pc)>0$, then $\theta(p')>0$ for some $p'<p$, and this
 contradiction implies Theorem 1. We make
 the required construction slightly more explicit as follows.

Recall the construction of Lemmas (18)--(21) of \BG.  We set $k=11$ and 
$\eta>0$; later we shall choose $\eta$ to be small.  Let $\varepsilon>0$ be
 such that $(1-\varepsilon)^{4k}>1-\eta$.  With this value of $\varepsilon$,
 we choose $r,L,K$ as in Lemma 4.1 and the preceding discussion, and we set
 $S=K+2r$.  Let ${\Cal R}^\pm = [-2L, 2L]\times V^\pm$ where
$$
V^\pm = \left\{(x_2, x_3) \in \Z^2 : \ 0 \leq  x_3 \leq (2k+2)S, -5L \pm
 \frac{L}{2S}x_3 \leq x_2 \leq 5L \pm \frac{L}{2S}x_3\right\}.
$$
We now define the {\it target zones\/} $V_{i,j}$ by
$$
V_{i,j} = w_{i,j} + [-L,L]\times [-2L, 2L] \times [0,2S]
$$
where $w_{i,j}=k(0,iL, 2jS)$ for $i,j \in \Z$ with $j \geq 0$ and $i+j$ even.
 Finally, let 
$$
{\Cal R} = \bigcup\Sb j \geq 0\\i+j \text{ even}\endSb \big\{({\Cal R}^+ \cup 
{\Cal R}^-)+w_{i,j}\big\}.
$$
These regions are illustrated in Figure 2.

\topinsert
\figure
\mletter{x_3}{5.5}{.25}
\lastletter{x_2}{9.8}{6.3}
\centerline{\epsfxsize=7cm\epsffile{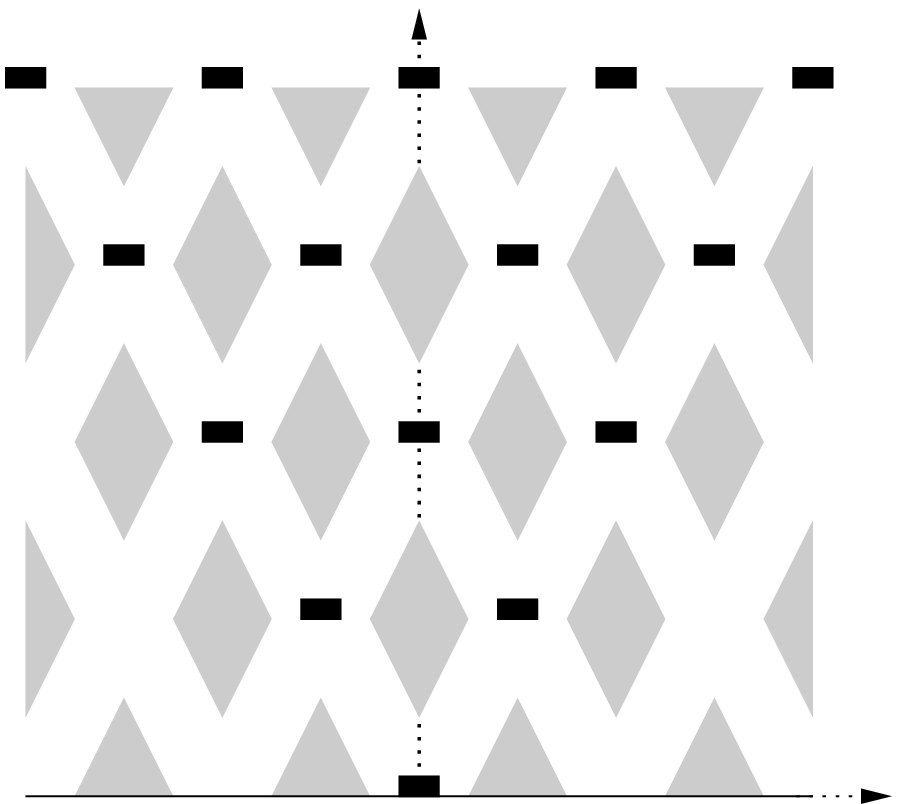}}

\capt{2}{The target zones are drawn in black,
and the region $\Cal R$ in white.}
\endinsert

The usual block variables are defined as follows.  Note that $D_r 
\subseteq V_{0,0}$.  We turn our attention to 
the target zones $V_{-1, 1}$ and $V_{1,1}$ 
and define indicator variables $\Xi_{-1,1}, \Xi_{1,1}$ by: $\Xi_{i,1}=1$ if
 and only if some vertex in $D_r$ is fully connected to some 
$r$-disk centred in $V_{i,1}$, by  paths contained 
entirely within the region ${\Cal R}$.  If $\Xi_{i,1}=1$, we let $\Delta_{i,1}$
 be an earliest $r$-disk centred in $V_{i,1}$ with the above property
 (`earliest' in order of third coordinate value).

We have by Lemma 4.1, the arguments of \BG, and the FKG inequality, that
$$
\P_p \left(\Xi_{-1,1} = \Xi_{1,1} = 1 \right) \geq (1-\varepsilon)^{4k} >
 1- \eta .
$$
The argument is now iterated from generation to generation.  Having
 constructed $\{\Xi_{i,j} : j \leq J\}$, we find the $\Xi_{i,J+1}$ by
 beginning with one of the  
$r$-disks $\Delta_{i-1,J}$,  $\Delta_{i+1,J}$  
already found in the construction yielding
 $\Xi_{i-1,J}$,  $\Xi_{i+1,J}$ (if both of these equal $0$, we set $\Xi_{i,J+1}=0$).
  If these variables both equal $1$, we choose the leftmost $r$-disk
 $\Delta$ of these two disks, 
and we declare $\Xi_{i,J+1}=1$ if and only if some point of $\Delta$
 is connected inside ${\Cal R}$ to every point of some $r$-disk
centred within the target zone $V_{i,J+1}$.
If $\Xi_{i,j}=1$, we say that the target zone $V_{i,j}$
has been {\it achieved\/} from the disk $D_r$.

This construction enables a comparison between the original process and a
 1-dependent directed conventional site percolation process on $\Z^2$ with
a certain intensity $1-\eta$.  If $\eta$ is sufficiently small, the latter process
 contains an infinite directed path with strictly positive probability and
 therefore so does the initial process.

We have shown that a supercritical directed percolation process
dominates, in a way made specific above, a 
{\it two-dimensional site\/} percolation
process with density close to 1. When proving Theorem 3 in Section 6,
we shall require the stronger property that it dominates
a certain {\it three-dimensional bond\/} 
percolation process. This will require some extra arguments.

We terminate this section with a note about the block construction
for the conventional percolation model on $\ldc$. The automorphism
group of $\ldc$ is rather different from that of $\ldca$, and this
has impact on the shape of the region corresponding to $B_{L,K}$ in
Lemma 4.1, and on the geometry of the ensuing block construction.
These turn out to be minor matters, and the results corresponding to
Theorem 1--3 are valid in this setting.

\section{5. Proof of Theorem 2} 
We begin with a subsidiary lemma.  For $x\in \Z^d$, we denote by $x(r)$ the
$r$-disk $x+D_r$.

\proclaim{Lemma 5.1}
Let $d\geq 2$, $p>\pc$, and $\xi>0$.  There exists a positive integer
 $R=R(p,\xi)$ such that
$$
\P_p \left(\exists z \text{\rm\ such that } x(r)\fc z(r),\  y(r)\fc z(r)\to 
\infty\right)>1-\xi
$$
for all $r \geq R$ and all distinct $x,y \in \Z^d$.
\endproclaim

\demo{Proof}
This may easily be shown when $d=2$, using path-intersection properties, and 
we therefore restrict ourselves to the case $d=3$ since this contains all the
 ingredients sufficient for the general case $d \geq 3$.

Let $\xi>0$, and find $\rho \in (0,1)$ such that a directed site
percolation process on $\vec\L^2$ with density $\rho$ is infinite with
probability at least $1-\half\xi$.  (An account
of the basic properties of directed percolation may be found in Section 1.6
of [\ci{G99}].) Set $k=11$
as in Section 4, and choose $\varepsilon>0$ such that:
$$
2k \varepsilon < \half \xi \tag 5.1
$$
and the block process referred to  at the end of Section 4 has density at
 least $1-\eta>\sqrt{\rho}$.  With this choice of $\varepsilon$, we choose $R$
 according to equation (4.1).  

Beginning with the disk $x(R)$, we make some initial steps of the block 
constructions of Section 4, but in a direction which carries us away from 
$y(R)$.  More specifically, we may assume without loss of generality that
$$
x_1-y_1 = \max \{\v x_i - y_i\v : 1 \leq i \leq 2\},\tag 5.2
$$
and for the moment we assume in addition that,
$$
\cases
\text{for all } u \in y + [-L,L] \times [-2L, 2L]\times [0,2S] 
 \text{ and}\\
 v \in x + (kL,0,2kS) + [-L,L] \times [-2L, 2L]\times [0,2S],\\
\lthreeca \text{ possesses a directed path from } u
 \text{ to } v.\endcases
\tag 5.3
$$ 
We shall see later how to adapt the proof when (5.3) fails.

As described in Section 4, with probability at least $(1-\varepsilon)^{2k
}>1-\half \xi$, $x(R)$  is connected to every point in some $R$-disk
centred in the region $x + (kL,0,2kS) + [-L,L] \times
[-2L,2L]\times [0,2S]$ by paths contained within the `tube' $x+ [-2L, 2L]
\times V^+$.  Let $A$ be the event that such an $R$-disk $w(R)$
can be
found, and we pick the earliest in the ordering induced by third coordinate
value.

Assume that the event $A$ occurs.  We now rotate our 
frames of reference and use the disks
$y(R)$ and $w(R)$ to initiate block constructions within disjoint subsets of
$\Z^3$, namely $y+\Cal R$ and $w+\Cal R$ where $\Cal R$
is given in Section 4.  
Since the difference in the first coordinates of $y$ and the centre
 of $w(R)$ exceeds $(k-1)L>9L$, and since the depth of ${\Cal R}$ is $4L$,
 this may be done.

Each step in the block constructions from $y(R)$ and $w(R)$ 
is successful with probability at least  $\sqrt{\rho}$,
whence {\it both\/} are successful with probability at least
$(\sqrt{\rho})^2=\rho$, which by assumption exceeds $\pc$.
It follows that the set of $(i,j)$ such
that both $V_{i,j}(y)$ and $V_{i,j}(w(R))$ are achieved from $y(R)$
and $w(R)$ respectively is infinite with probability at least
$1-\half \xi$.

We call the pair $(i,j)$ {\it green\/} if, for all $u\in V_{i,j}(y)$
and $v\in V_{i,j}(w(R))$, $u\to v$ 
in the convex hull of $V_{i,j}(y)$ and 
$V_{i,j}(w(R))$.  We have by assumption (5.3) that
$$
\gamma = \P_p \big((i,j)\text{ is green}\big)>0,
$$
and it follows by an application of the FKG inequality that, with probability
(conditional on $A$) at least $1-\half \xi$, there exists a 
green $(i,j)$ such that: $V_{i,j}(y)$
 is achieved from $y(R)$ and $V_{i,j}(w(R))$ is achieved from $w(R)$,
and in addition the block constructions from $y(R)$
and $w(R)$ include infinite paths of blocks beginning respectively at
$V_{i,j}(y)$ and $V_{i,j}(w(R))$.

If $A$ occurs, and also the last event, 
there exists $z$ such that $x(R)\fc z(R)\to\infty$ 
and $y(R)\fc z(R)$.  The probability of failure does not exceed 
$1-\P_p(A) + \half \xi<\xi$, and the claim is proved.

Finally we return to assumption (5.3).  When (5.3) fails, we need to continue
 the construction of the event $A$ by adding further steps to the block
 construction from $x(R)$ until we obtain a target zone $V_{I,J}(x)$ with the
 property that, for all $u\in V_{0,0}(y)$, $v\in V_{I,J}(x)$ there exists a
 directed path of $\lthreeca$ from $u$ to $v$.  We let $A$ be the
 event that the construction successfully attains an $r$-disk
centred in $V_{I,J}(x)$ and then we argue as before.  In this case,
 $\P_p(A) \geq (1-\varepsilon)^{2kJ}$, and we amend the choice of
 $\varepsilon$ accordingly.
\qed\enddemo

\demo{Proof of Theorem 2}
(a) Let $p_c <\alpha \leq \beta \leq 1$, $\xi>0$, and pick $R=R(\alpha,\xi)$
 according to Lemma 5.1.  Let $\Delta^\gamma_u$ denote the earliest $R$-disk
such that the vertex
$u$ is fully $\gamma$-connected to $\Delta_u^\gamma$.
Here, `earliest' means in the ordering induced by the
third coordinate value, and, if there is a choice, we take the earliest in
some predetermined ordering.  Let $A^\gamma_u$ be
the event that such a $\Delta^\gamma_u$ 
exists.
It is elementary, by a `re-start' argument, that
$$
\P (A^\gamma_x \mid x \overset \gamma \to \to\infty)=1 \quad \text{ for }
 \gamma>p_c.
$$ 
By Lemma 5.1,
$$
\P\big(\exists z\text{ such that } x\overset \alpha\to\to z\overset 
\alpha\to\to\infty,\ 
y\overset \beta\to\to z \mid A^\alpha_x \cap A^\beta_y\big) > 1 - \xi,
$$
whence
$$\align
\P\big(\exists z\text{ such that } x\overset \alpha\to\to z
  \overset \alpha\to\to\infty,
 \ y\overset \beta\to\to z\big) & >  (1-\xi)\P (A^\alpha_x \cap
 A_y^\beta)\\
& \geq (1-\xi) \P(x\overset\alpha\to\to \infty, \ y 
\overset \beta\to\to\infty).
\endalign
$$
Since this holds for all $\xi >0$, the claim of part (a) follows.

\flushpar (b) The right-continuity of $\theta$ is an immediate 
consequence of the fact 
that $\theta$ is a decreasing limit of continuous non-decreasing functions
 (cf.\ [\ci{G99}], p.\ 203).  Since $\theta (\pc)=0$, by Theorem 1, it follows that
 $\theta$ is continuous on $[0,\pc]$.  In order to prove the left-continuity 
of $\theta$ on
the interval $(\pc,1]$, one adapts the 
argument of  [\ci{BeKea}] in the usual way (see 
[\ci{G99}], p. 203), making use of the result of part (a).
\qed\enddemo

\section{6. Proof of Theorem 3}
Rather than developing the argument of [\ci{GKZ}], we make use of a result
of [\ci{BenPP97}], where it is proved in the context of the conventional
model on $\vec\L^3$ that, for sufficiently large $p$, the undirected
graph $C$ is transient almost surely on the event $\{|C|=\oo\}$.
The further question is posed in [\ci{BenPP97}] whether this conclusion
is valid under the weaker hypothesis that $p$ exceeds the appropriate critical
probability, and our Theorem 3 is a positive answer to this question. Although
Theorem 3 as stated relates to the alternative model, similar arguments apply
for the conventional model.

As explained earlier, we do not include all the details of the required proof,
but instead we describe only the salient features. We apologise to
those who might have savoured the complicated notation 
of an overlong full proof, but we hope that readers familiar
with [\ci{BG1}] will agree with our decision.

We shall construct a block process whose target zones are indexed 
in the following manner.
Let $G$ be the graph having vertex set $\Z^3$ and edge set as follows: for
$x=(x_1,x_2,x_3)$ and $y=(y_1,y_2,y_3)$, we place a directed edge from
$x$ to $y$ if and only if $|x_1-y_1|+|x_2-y_2|=1$ and $y_3-x_3=1$.
Note that $G$ is a subgraph of $\lthreeca$. The target zones
of our block process will be indexed by the set $W$ of vertices of $G$
which are accessible along directed paths from the origin, and the block
process will proceed by building connections within regions of $\lthreeca$
represented by the edges joining vertices in $W$.

Suppose that $p>\pc$. Let $\eps>0$, and choose $r$, $L$, $K$ as
in (4.1) and Lemma 4.1. We shall assume a bound on $\eps$ of the 
form $\eps<\eps_0$,
where the small quantity $\eps_0$ will be chosen later.
By the usual `re-start argument', if $0\to\oo$,
there exists almost surely an $r$-disk $\Delta$ such that 
$0$ is fully connected to $\Delta$. We use the argument of Section 4
to construct a block process from $\Delta$ in which the target zones are indexed
by $W$. This block
construction dominates (stochastically) a directed
site percolation model on $W$ with some density $1-\eta(\eps)$, where
$\eta(\eps)\downarrow 0$ as $\eps \downarrow 0$. Assume for the moment
that we may replace the word `site' with the word `bond' in
the last sentence. It may then be shown, by Proposition 1.2 and 
Theorem 1.3 of [\ci{BenPP97}] together with some
standard arguments concerning electrical networks (see [\ci{DoyS}, \ci{GKZ}]), 
that $C$ is transient with (conditional) 
probability
at least $1-\xi(\eps)$, for some $\xi(\eps)$ satisfying 
$\xi(\eps)\downarrow 0$ as $\eps\downarrow 0$.
Since $\eps$ was arbitrary, the conclusion will follow. There are two gaps
in this argument, namely:
\roster
\item"{(a)}" to show that supercritical directed percolation
on $\lthreeca$ dominates a certain directed {\it bond\/} 
percolation process on $W$ having large density,
\item"{(b)}" to check that the conclusion of Theorem 1.3 of [\ci{BenPP97}] 
is valid for the graph $G$.
\endroster

An examination of the proof of
Theorem 1.3 of [\ci{BenPP97}] reveals that it is easily adapted to the graph
$G$,
and one may conclude that there exists a probability measure on directed
paths from the origin in $G$ that has exponential
intersection tails. We turn therefore to point (a).

We begin with a `local connection lemma'. Let $p$, $\eps$ ($<\eps_0$),
$r$, $L$, $K$ be as above, and let $S=K+2r$. Let $B_{l,k}=[-l,l]^2\times[0,k]$
as before,
and let $\Cal D$ be the set of all $r$-disks centred in $B_{3L,3K}$.
For $M\ge 5\max\{K,L\}$ and $\Delta, \Delta'\in\Cal D$, let $E_{\Delta,\Delta'}$
be the event that there exists an $r$-disk $\Delta''$ 
contained entirely in $B_{M,2M}$ such
that both $\Delta$ and $\Delta'$ are fully connected to 
$\Delta''$ in $B_{M,2M}$.
This event is illustrated in Figure 3.

\topinsert
\figure
\mletter{B_{M,2M}}{7.6}{2}
\mletter{B_{3L,3K}}{4.7}{3.65}
\mletter{\Delta''}{7.2}{2.7}
\mletter{\Delta}{5.6}{5.5}   
\lastletter{\Delta'}{6.3}{5.07}
\centerline{\epsfxsize=7cm\epsffile{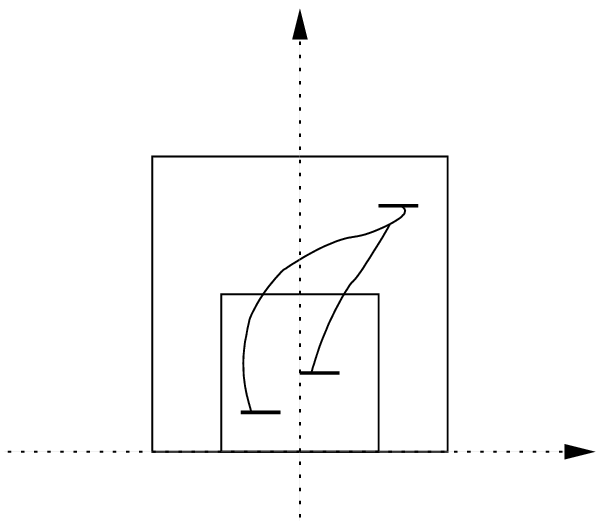}}

\capt{3}{An illustration of the event $E_{\Delta,\Delta'}$.}
\endinsert

\proclaim{Lemma 6.1} There exist an integer $M(\eps)\ge 5\max\{K,L\}$ and a 
function $\b(\eps)$ satisfying
$\b(\eps)\downarrow 0$ as $\eps\downarrow 0$ such that
$$
\P_p(E_{\Delta,\Delta'})\ge 1-\b(\eps)\q\text{for all } \Delta,\Delta'\in\Cal D,
$$
whenever $M\ge M(\eps)$.
\endproclaim

\demo{Proof} 
Let $x(r)$ denote the $r$-disk centred at the vertex $x$. For $x,y\in\Z^3$,
let $F_{x,y}(N)$ be the event that there exists an $r$-disk $\Delta''$
such that $x(r)\fc\Delta''$ and $y(r)\fc \Delta''$ in $x+B_{N,N}$.
By Lemma 5.1, there exists $\gamma(\eps)$ satisfying $\gamma(\eps)
\downarrow 0$ as $\eps\downarrow 0$ such that
$$
\P_p\Bigl(\lim_{N\to\oo} F_{x,y}(N)\Bigr)\ge 1-\g(\eps)\q\text{for
all } x,y\in\Z^3.
$$
Therefore there exists $N=N_{x(r),y(r)}$ such that
$$
\P_p\bigl( F_{x,y}(N_{x(r),y(r)})\bigr)\ge 1-2\g(\eps)\q\text{for
all } x,y\in\Z^3.
$$
Let $\Delta,\Delta'\in\Cal D$, take $x(r)=\Delta$,
$y(r)=\Delta'$, and choose $M=M(\eps)\ge 5\max\{K,L\}$ sufficiently
large that $\Delta+B_{N_{\Delta,\Delta'},N_{\Delta,\Delta'}}
\subseteq B_{M,2M}$ for all
$\Delta\in\Cal D$. With this value of $M$,
$\P_p(E_{\Delta,\Delta'})\ge 1-2\gamma(\eps)$
uniformly in $\Delta$ and $\Delta'$.
\qed\enddemo

We illustrate next how Lemma 6.1 may be used to 
show that supercritical directed percolation
on $\lthreeca$ dominates a certain {\it two-dimensional\/}
directed bond percolation process. It is easier to draw pictures in 
this case, and it will be explained later how to extend the claim to 
the three-dimensional graph $G$.

Let $1-\rho(\eps)$ be the probability that the 
two-dimensional block construction
of Section 4, initiated from the disk $D_r$, yields
an infinite structure. Since this block process dominates a
directed
site percolation process having some density $1-\eta(\eps)$
where $\eta(\eps)\downarrow0$ as $\eps\downarrow 0$, we have by 
standard arguments 
(see, for example,  [\ci{Dur84}, \ci{Lig}] and the references in
[\ci{BG1}])
that there exist $v(\eps),\gamma(\eps)\ge 0$ such that:
\roster
\item"{(a)}" $v$ is non-increasing in $\eps$, and strictly positive
when $\rho(\eps)<1$,
\item"(b)" $\gamma(\eps)\downarrow 0$ as $\eps\downarrow 0$,
\item"(c)" for $i,j\in \Z$ such that $j\ge 0$ and $i+j$ is even,
if $|i|/j< v(\eps)$ then the target zone $V_{i,j}$ is
achieved with probability at least $1-\gamma(\eps)$ by an open path
lying entirely within the convex region 
$$
\Cal C(\eps)= B_{5L,10L}+
\left\{(0,x_2,x_3): |x_2|/x_3\le \frac{L}{2S}\cdot 2v(\eps)\right\}.
$$
\endroster
Let $M=M(\eps)$ and $\b(\eps)$ be given as in Lemma 6.1. It follows that
we may find positive integers $I$, $J$ depending
on $\eps$ such that:
\roster
\item"(i)" $IL/M$ and $JS/M$ are large, say $IL/M,JS/M\ge 10$,
\item"(ii)" $v(\eps)<I/J<2v(\eps)$,
\item"(iii)" for all $r$-disks $\Delta$ lying in $B_{M,2M}$,
there exists with probability at least 
$1-\g(\eps)$ an $r$-disk $\Delta'$ 
centred in $11(0,IL,2JS)+B_{3L,3K}$ such that $\Delta\fc\Delta'$
in $\Cal E(\eps)=B_{M,2M}+\Cal C(\eps)$.
\endroster

\midinsert
\figure   
\mletter{x_3}{6.4}{.3}
\mletter{z_{-2,2}}{3}{2.5}
\mletter{z_{0,2}}{6.3}{2.5}
\mletter{z_{2,2}}{9}{2.5}
\mletter{Z_{-1,1}}{3.2}{3.65}
\mletter{Z_{1,1}}{8.5}{3.65}
\mletter{z_{-1,1}}{4.4}{4.65}
\mletter{z_{1,1}}{7.38}{4.65}
\mletter{Z_{0,0}}{6.9}{5.75}
\mletter{x_2}{10.3}{6.35}
\lastletter{0}{6.32}{6.95}
\centerline{\epsfxsize=9cm\epsffile{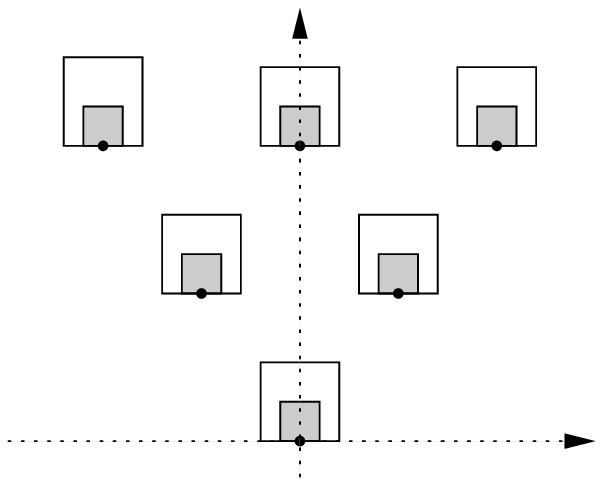}}

\capt{4}{The target zones $Z_{i,j}$. Each smaller box $Z_{i,j}$
is a translate of
$B_{3L,3K}$, and is contained in a translate of $B_{M,2M}$.}
\endinsert

For integers $i,j$ such that $j\ge 0$ and $i+j$ is even,
we define  target zones
$$
Z_{i,j}=z_{i,j}+B_{3L,3K},
$$
where $z_{i,j}=11(0,iIL,2jJS)$. See Figure 4.

We now define new block indicator variables $\Theta_{i,j}$ inductively
as follows. We set $\Theta_{0,0}=1$.
For $i=\pm 1$, we set $\Theta_{i,1}=1$ if and only if
$D_r$ is fully connected to some $r$-disk centred in
$Z_{i,1}$ by open paths in $\Cal E(\eps)$. 
If this holds, we let $\Delta_{i,1}$ be the earliest
such $r$-disk in $Z_{i,1}$. Having constructed $\{\Theta_{i,j}: j\le R\}$,
we find the $\Theta_{i,R+1}$ as follows. Let $g(i,R+1)$ be the set
of all $i'\in\{i-1,i+1\}$ such that
$\Theta_{i',R}=1$. If $g(i,R+1)$ is empty, we set $\Theta_{i,R+1}=0$.
If $g(i,R+1)$ contains a singleton, say the value $i'$, we
set $\Theta_{i,R+1}=1$ if and only if $\Delta_{i',R}$
is fully connected to some $r$-disk centred
in $Z_{i,R+1}$ by open paths
lying in $z_{i',R}+\Cal E(\eps)$, and we denote by  $\Delta_{i,R+1}$
the earliest such $r$-disk. So far we have retained much of the manner of
the construction given in Section 4. 
However, an important  difference arises when $g(i,R+1)$
contains {\it both\/} $i-1$ and $i+1$. In this case we set $\Theta_{i,R+1}
=1$ if and only if the following occur:
\roster
\item"1." for $i'= i-1$ and $i'=1+1$, the $r$-disk $\Delta_{i',R}$ is
fully connected to some $r$-disk centred in
$Z_{i,R+1}$ by open paths lying in $z_{i',R}+\Cal E(\eps)$,  
\item"2." writing $\Delta(i')$ for the earliest such $r$-disk referred
to above, the event $E_{\Delta(i-1),\Delta(i+1)}$
occurs.
\endroster
If these occur, we denote by $\Delta_{i,R+1}$ the earliest
$R$-disk attained in the definition of $E_{\Delta(i-1),\Delta(i+1)}$.
This event is illustrated in Figure 5.

\midinsert
\figure
\mletter{\Delta_{i,R+1}}{5.8}{.6}
\mletter{Z_{i-1,R}}{1.6}{5.7}
\mletter{Z_{i,R+1}}{4}{2}
\mletter{\Delta_{i-1,R}}{3.4}{6.6}
\mletter{Z_{i+1,R}}{6.25}{5.7}
\lastletter{\Delta_{i+1,R}}{8.1}{6.8}
\centerline{\epsfxsize=7cm\epsffile{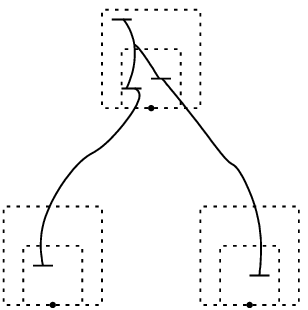}}

\capt{5}{An illustration of the definition of the block variable
$\Theta_{i,R+1}$.}
\endinsert

The $\Theta_{i,j}$ are dependent random variables, but the extent of
their dependence is limited. Suppose
$\{\Theta_{i,j}:j\le R\}$ have been observed and the corresponding disks
$\Delta_{i,j}$ found. The pairs $(\Theta_{i,R+1}, \Delta_{i,R+1})$,
$-R-1\le i \le R+1$, have some interdependence 
owing to the fact that the open paths from different
$\Delta_{i,R}$ may lie close to one another.
Using conditions
(i)--(iii) together with the fact that $v(\eps)$ is decreasing in $\eps$,
there exists
an integer $T=T(\eps_0)$ such that: conditional on
the pairs $(\Theta_{i,R},\Delta_{i,R})$, $-R\le i\le R$,
the family of pairs $(\Theta_{i,R+1},\Delta_{i,R+1})$, $-R-1\le i\le R+1$,
is $T$-dependent (see [\ci{G99}], page 178 for a definition
of $T$-dependence). This observation is illustrated in
Figure 6.

\topinsert
\figure
\centerline{\epsfxsize=9cm\epsffile{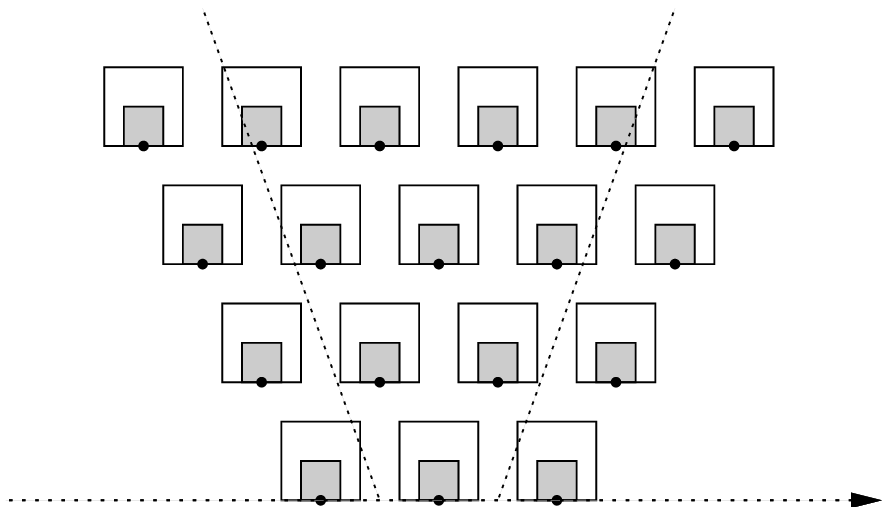}}

\capt{6}{Since we restrict ourselves to open paths lying within a certain
`wedge', the dependence between block variables
has a range which is bounded in $\eps$ ($<\eps_0$).}
\endinsert

Each step in the above
construction of the 
$\Theta_{i,j}$ is successful with probability at least
$1-2\gamma(\eps)-\b(\eps)$, 
which approaches $1$ as
$\eps\downarrow 0$. Since $T$ is an absolute constant,
we deduce by the comparison theorem
of [\ci{LSS}] (see also Theorem (7.65)
of [\ci{G99}]), that the $\{\Theta_{i,j}\}$, together with
the successful connections between the $r$-disks $\Delta_{i,j}$,
dominate
(stochastically) the open cluster at the 
origin of a directed {\it bond\/} percolation
process on $\vec\L^2$ having density approaching 1 as
$\eps\downarrow 0$.

Several details are missing from the foregoing argument, 
of which one is an account
of the `steering' necessary to achieve property (iii) above. This
follows a standard route, and is omitted.

We return finally to point (a) before Lemma 6.1.
It is required to show that the two-dimensional 
construction of the $\Theta_{i,j}$ may be
extended to a three-dimensional construction with target 
zones indexed by the set $W$. This we achieve with the aid of some pictures.
The target zones are now
$$
Z_{i,j,k}=z_{i,j,k}+B_{3L,3K},
$$
where $z_{i,j,k}=11(iIL,jIL,4kJS)$, as
$(i,j,k)$ ranges over integer vectors satisfying
$|i|\le k$, $|j|\le k$, and $i+k$ and
$j+k$ are even.
A plan of these target zones is drawn in Figure 7.
Boxes of the form $11(iIL,jIL,2kJS)+B_{3L,3K}$ for odd values of $k$ are
termed `intermediate zones'.

\topinsert
\figure
\mletter{x_2}{1.75}{1.5}
\mletter{x_1}{3.1}{2,6}
\mletter{k=0}{1.8}{4.6}
\mletter{k=1}{4.75}{4.6}
\lastletter{k=2}{8.5}{4.6}
\centerline{\epsfxsize=10cm\epsffile{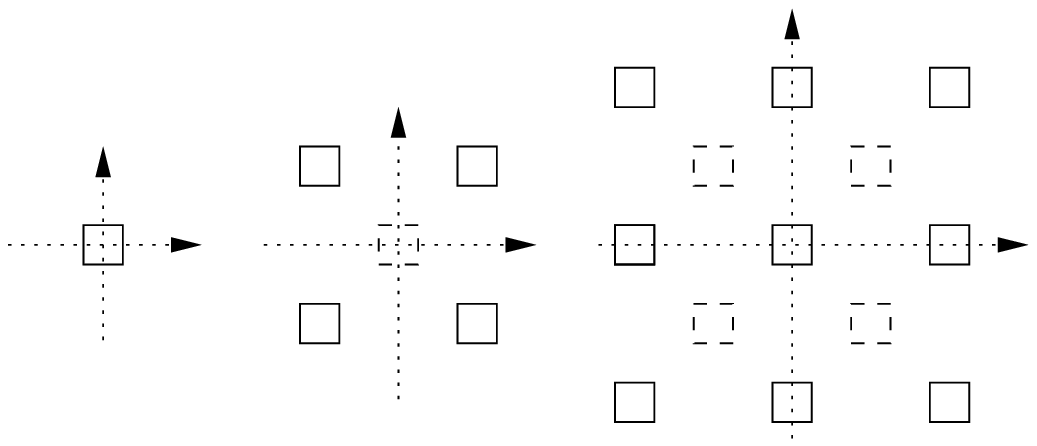}}

\capt{7}{A plan of the target zones $Z_{i,j,k}$.}
\endinsert

We describe next the open connections sought in defining the
block indicator variables $\Theta_{i,j,k}$. First we set
$\Theta_{0,0,0}=1$. For $i,j\in\{-1,1\}$,
we declare $\Theta_{i,j,1}=1$ if and only if:
\roster
\item"(a)" there exists an $r$-disk $\Delta$ centred in the intermediate zone
$11(iL,0,2JS)+B_{3L,3K}$ such that $D_r\fc \Delta$ in a certain
convex region $\Cal E$, and
\item"(b)" there exists an $r$-disk $\Delta'$ in $Z_{i,j,1}$ such that
$\Delta\fc\Delta'$ within a certain translate/rotation of $\Cal E$.
\endroster
The region $\Cal E$ corresponds to the $\Cal E(\eps)$ used above.
When such connections exist, we let $\Delta_{i,j,1}$ be the
earliest disk in $Z_{i,j,1}$ which is thus
reached from $D_r$.

The next step is similar to that described in the two-dimensional case.
We omit most of the details, but concentrate on one
illustrative example. Suppose for the sake of
the illustration that $\Theta_{i,j,1}=1$ for all $i,j$, and that
we are seeking a definition of $\Theta_{0,0,2}$. We set $\Theta_{0,0,2}=1$
if and only if the following hold:
\roster
\item"(i)" there exists an $r$-disk $\Delta(-1,1)$ (respectively $\Delta(1,1)$)
centred in the intermediate zone $11(0,IL,6JS)+B_{3L,3K}$
such that $\Delta_{-1,1,1}\fc \Delta(-1,1)$ (respectively $\Delta_{1,1,1}
\fc\Delta(1,1)$) by paths lying within a certain translate/rotation
of $\Cal E$,
\item"(ii)" there exists an $r$-disk $\Delta_1$ 
such that $\Delta(-1,1)\fc \Delta_1$ and $\Delta(1,1)
\fc\Delta_1$ in $11(0,IL,6JS)+B_{M,2M}$,
\item"(iii)" there exists an $r$-disk $\Delta(-1,-1)$ (respectively $\Delta(1,-1)$)
centred in the intermediate zone $11(0,-IL,6JS)+B_{3L,3K}$
such that $\Delta_{-1,-1,1}\fc \Delta(-1,-1)$ (respectively $\Delta_{1,-1,1}
\fc\Delta(1,-1)$)  
by paths lying within a certain translate/rotation
of $\Cal E$,
\item"(iv)" there exists an $r$-disk $\Delta_{-1}$ 
such that $\Delta(-1,-1)\fc \Delta_{-1}$ and $\Delta(1,-1)
\fc\Delta_{-1}$ in $11(0,-IL,6JS)+B_{M,2M}$,
\item"(v)" there exists an $r$-disk $\Delta_1'$ (respectively
$\Delta_{-1}'$) centred in the box $z_{0,0,2}+B_{3L,3K}$
such that $\Delta_1\fc\Delta_1'$ (respectively
$\Delta_{-1}\fc\Delta_{-1}'$) by paths lying within certain translates/rotations
of $\Cal E$,
\item"(vi)" there exists an $r$-disk $\Delta$ in $Z_{0,0,2}$
such that $\Delta_1'\fc\Delta$ and $\Delta_{-1}'\fc\Delta$ in $Z_{0,0,2}$.
\endroster
If these events occur, we define $\Theta_{0,0,2}=1$,
and we let $\Delta_{0,0,2}$ be the earliest $r$-disk $\Delta$ thus
accessible in (v).

One may define the variables $\Theta_{i,j,k}$ in a similar inductive way.
As before, it may be shown via the stochastic domination
theorem of [\ci{LSS}] that the block process dominates
a bond percolation process with high density,
and our sketch of the proof
of Theorem 3 is complete.

\section{7. Infinite-range percolation}
We consider next a long-range directed percolation model which extends our 
earlier results for the alternative lattice ${\ldca}$. Let 
$d\ge 2$, and let $\bp = \{p(\bx)$:
 $\bx \in \Z^{d-1}\}$ be a collection of numbers satisfying $0 \leq p (\bx)<1$.
  Consider the vertex set $\Z^d$, and write $x=(x_1,x_2,\dots,x_d)
\in \Z^d$ as $x=(\bx,t)$ where
 $\bx = (x_1, x_2,\dots, x_{d-1})\in \Z^{d-1}$ and $t=x_d\in \Z$.  For every
 $\bx, \boldy \in \Z^{d-1}$ and $t\in \Z$, we place a directed edge from $(\bx,t)$ 
to $(\boldy, t+1)$ with probability $p(\boldy-\bx)$.  Each such pair is
 joined by an edge independently of the presence or absence of other edges.

We shall use the same notation as earlier.  For example, we write 
$x\to y$ if there exists a directed path from $x$ to $y$, and we let 
$\theta (\bp) = \P_{\bp}(0\to\infty)$, where $\P_{\bp}$ 
denotes the appropriate probability measure.

We shall require a certain amount of symmetry, and to this end we shall
 assume that:
\roster
\item"{(a)}" $\bp$ is invariant under sign changes of components, in that
 $p(\bx)=p(\bx')$ whenever $\bx'$ is obtained from $\bx$ by changing the signs
 of any of the $d-1$ components of $\bx$, and

\item"{(b)}" $\bp$ is invariant under permutations of components, in that 
$p(\bx) = p(\pi\bx)$ where $\pi$ is a permutation of $1,2,\dots, d-1$, and
$\pi(x_1, x_2,\dots, x_{d-1}) = (x_{\pi(1)}, x_{\pi(2)}, \dots, x_{\pi(d-1)})$.
\endroster 

\flushpar The {\it range\/} $R$ of the process is 
defined as $R=\sup\{\v\bx\v:p(\bx)>0\}$,
 and the  process  is said to have {\it infinite range\/} if $R=\infty$.

Under what further assumptions on $\bp$ may one adapt the ideas underlying the
 block construction of Section 4?  We will present sufficient conditions on $\bp$,
 and it will follow in the usual way that, when these conditions are valid,
 the usual conclusions follow, including the continuity of slab critical
 points, the (suitably generalized)
uniqueness result of Theorem 2(a), the fact that the {\it critical\/}
 process (for a suitable parametric family of functions $\bp$
satsfying a further condition of continuity) dies out, and
 many other observations.  The details of such consequences are not included
 here, since they follow already familiar lines.  Instead, we make specific
 our sufficient conditions on $\bp$, and we outline the steps which follow for the
 required block construction.

Infinite-range {\it undirected\/} percolation has been studied in one and 
higher dimensions in [\ci{AN86}, \ci{GKM}], and a block
construction was developed in [\ci{MeeS}] subject to a rather severe 
condition on the decay rate of probabilities of long-range edges.  Infinite-range
{\it directed\/} percolation, and particularly some of the facts referred to
above, are relevant to the scaling limit of critical directed percolation in
high dimensions, proved in [\ci{vdHS00}, \ci{Sl00}]. The decay
rate required by the forthcoming conditions on
the function $\bp$ is substantially weaker than that
assumed in [\ci{MeeS}] for undirected percolation.

For $\bx=(x_1,x_2,\dots,x_{d-1})\in\Z^{d-1}$, we define
$\v\bx\v_\infty = \max \{\v x_i\v:1 \leq i \leq d-1\}$,
and we write
$$
\Sigma_u=\sum_{\bx:|\bx|_\oo>u}p(\bx)\q\text{for } u>0.
$$
The relevant
conditions on $\bp$ are the following.
\roster
\item"{I.}" {\it Summability.} $\sum_{\bx\in \Z^{d-1}} p(\bx)<\infty$.

\item"{II.}" {\it  Aperiodicity.} For every $\bx \in \Z^{d-1}$, the
 set $\{t:\P_{\bp} (0\to (\bx,t))>0\}$ has greatest common divisor 1.

\item"{III.}"  {\it Tail regularity.} There exists an integer $\alpha>1$
 and a real $\xi\in(0,1)$
 such that $\Sigma_{\alpha h} \le \xi\Sigma_h$
for all $h\in\R$ satisfying $h\ge 1$.
\endroster

We next discuss these conditions.  Condition I holds if and only if every
 vertex has almost surely finite vertex degree.  Condition II is a convenience
 but is not essential.  If Condition II fails, 
then one would sometimes need to restrict
 oneself to an appropriate subset of $\Z^d$.  Condition III is a condition of
smoothness on the manner in which 
$p(\bx)$ decays for large $\v\bx\v_\oo$, and this will be
 required for the renormalization argument. 

Condition III is not over severe. Assume for the sake
of illustration that
$p(\bx)=g(|\bx|_\oo)$ for large $|\bx|_\oo$, where $g(v)=v^{-\beta}$.
Then Condition III is satisfied whenever Condition I holds, namely if 
$\b>d-1$. The condition is however not satisfied if, for example,
$g(v)=v^{-d+1}(\log v)^{-2}$.  

It is easy to see, 
by the symmetry
of $\bold p$, that Condition III implies
$$
\sum\Sb \bx:|\bx|_{\infty}\leq \alpha h \\ x_1>h\endSb p(\bx) \geq
\frac{\Sigma_h-\Sigma_{\a h}}{2(d-1)} \ge 
\frac{(1-\xi)\Sigma_h}{2(d-1)}\ge  \frac{1-\xi}{2(d-1)} \sum_{\bx:x_1>h}p(\bx),
\tag 7.1
$$
and we shall make use of this fact later.

Let us assume henceforth that $\bp$ satisfies 
Conditions I, II, III for some pair $\alpha,
\xi$, and we assume as before that $d=3$.  We claim that the `usual
 theorems' follow, and to this end we now sketch the necessary extra steps in
 order to achieve a block construction as in Section 4.  The principal step 
is to  establish an equivalent of Lemma 4.1.  Whereas Lemma 4.1 concerned the
 numbers of points on the top and sides of a block  $B_{L,K}$ which are
 endpoints of directed paths of the block originating
from the disk $D_r$, in
 the infinite-range setting we concentrate on the number $N_\T(L,K)$ of 
such points
 on the top of $B_{L,K}$, and the mean number 
$\widetilde{N}_\S(L,K)$ of edges exiting $B_{L,K}$ by
 its sides.  That is, we replace $N_\S (L,K)$ by
$$
\widetilde{N}_\S(L,K) = \sum\Sb x=(\bx,t)\\ y=(\boldy, t+1)\\0\le t<K\endSb 
p(\boldy-\bx)\tag 7.2
$$
where the summation is over $x \in B_{L,K}$ such that $D_r\to
 x$ in $B_{L,K}$, and over $y\in \{[-\gamma L, \gamma L]^2\setminus
 [-L,L]^2\}\times [1,K]$, where $\gamma=2\alpha-1$.  
For $\bv \in \{-1,+1\}^3$, we consider the
 sub-facet of the sides of $B_{L,K}$ indexed by $\bv$, and we define
 $\widetilde{N}^{\bv}_\S(L,K)$ as follows. For simplicity, assume this
 sub-facet is $[0,L]\times \{L\} \times [0,K]$, and let
$$
\widetilde{N}_\S^{\bv} (L,K) = \sum\Sb x=(\bx,t)\\y=(\boldy,t+1)\\0\le t<K\endSb
\smash{\hskip-4mm{}^{\bv}}\ p(\boldy-\bx)
\tag 7.3
$$
where the summation is as in (7.2) but with $y$ restricted to the region
 $[0,\gamma L]\times (L,\gamma L]\times [0,K]$. See Figure 8.

\topinsert
\figure
\mletter{x_2}{5.7}{.5}
\mletter{(\gamma L,\gamma L)}{9.2}{.8}
\mletter{x_1}{9.4}{4}
\mletter{\bx}{5.2}{3.4}
\mletter{\boldy}{8.4}{1.23}
\lastletter{(L,-L)}{7.6}{5.8}
\centerline{\epsfxsize=7cm\epsffile{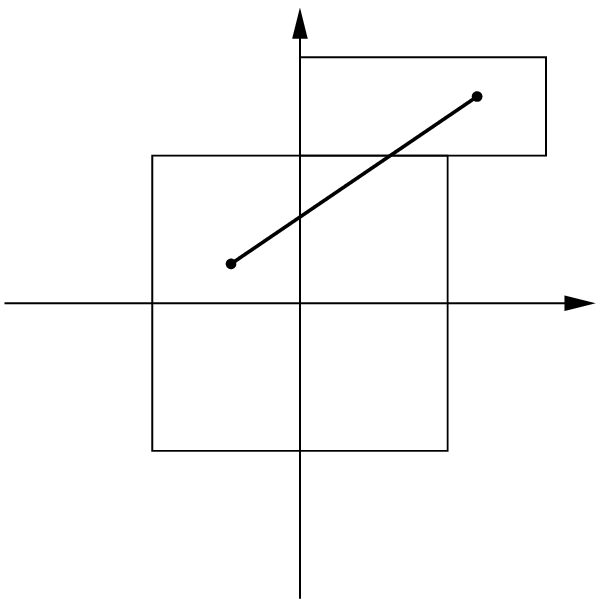}}

\capt{8}{A map of $(x_1,x_2)$-space. The set of possible
values of $\boldy$ in (7.3) is in the upper rectangle.}
\endinsert

With a few minor changes, the proof of Lemma 4.1 may be adapted to obtain
 the corresponding conclusion with $N^\bv_\S$ replaced by
 $\widetilde{N}^\bv_\S$.  We highlight one difference, namely that related to
Condition III.  
For $\bx\in[-L,L]^2$, we set
$$
\nu_L(\bx)=\sum_{\boldy\notin[-L,L]^2} p(\boldy-\bx),
\tag 7.4
$$
and
$$
Q=\sum\Sb x=(\bx, t)\\0\le t<K\endSb \nu_L(\bx)
=\sum\Sb x=(\bx, t)\\ y=(\boldy, t+1)\\0\le t<K\endSb p(\boldy - \bx)
\tag 7.5
$$
where the sums are over $x\in B_{L,K}$ such that $D_r\to x$ in
 $B_{L,K}$, and over $y\in \{(-\infty, \infty)^2\setminus[-L,L]^2\}
\times [1,K]$. We require a statement of the form: $\wt N_\S(L,K)$ is large
if and only if $Q$ is large.

\proclaim{Lemma 7.1} 
We have that $\frac18(1-\xi) Q\le \wt N_\S(L,K)\le Q$.
\endproclaim
 
\demo{Proof}
For $\bx\in[-L,L]^2$, let
$$
\nu_L^\g(\bx)=\sum_{\boldy\in[-\g L,\g L]^2\setminus[-L,L]^2} p(\boldy-\bx).
\tag 7.6
$$
For $i=1,2$ and $\eta=\pm$, let $\nu_L^{i,\eta}(\bx)$ (respectively,
$\nu_L^{\g,i,\eta}(\bx)$) be given as in
(7.4) (respectively, (7.6)) but with the restriction 
to vertices $y$ satisfying $\eta y_i> L$.
Evidently, by (7.1),
$$
\nu_L^\g(\bx) \ge \tfrac12\sum_{i,\eta} \nu_L^{\g,i,\eta}(\bx) 
\ge  \tfrac18(1- \xi)\sum_{i,\eta} \nu_L^{i,\eta}(\bx)
\ge \tfrac18 (1-\xi)\nu_L(\bx).
$$
Therefore,
$$
\wt N_\S(L,K) =\sum\Sb x=(\bx, t)\\0\le t<K\endSb  \nu_L^\g(\bx)
\ge \sum\Sb x=(\bx, t)\\0\le t<K\endSb \tfrac18(1-\xi)\nu_L(\bx)
=\tfrac18(1- \xi) Q,
$$
as required, where the sums are over $x\in B_{L,K}$ such that $D_r\to x$ in
 $B_{L,K}$.
The other inequality of the lemma is a triviality.
\qed\enddemo

Having achieved an infinite-range version of Lemma 4.1, one embarks on the 
block construction, and it is here that one sees the value of the definition
 of the $\widetilde{N}^{\bv}_\S(L,K)$.  The maximum lateral extension of the
 points $y$ in (7.2) is $\g$ times the dimension of the sides of
 $B_{L,K}$, where $\g$ is a constant which is independent of $L$ and $K$.
  It is the last fact which enables a block construction in which each step
 uses a number $l$ of extensions of the basic step of Lemma 4.1, the number
 $l$ being independent of $L$ and $K$.  With some care but with no further
 ideas of consequence, it follows that, 
with an appropriate choice of parameters,
 the block process stochastically dominates 
a conventional directed site percolation
 process having density close to $1$.

We finish with a discussion of the infinite-range analogue of the statement
 $\theta (\pc)=0$.  Let $\{\bp^r:0 \leq r \leq 1\}$ be a parametric family of
 functions satisfying conditions (a) and (b) at the beginning of this section,
 and suppose that every $\bp^r$ satisfies conditions I, II, III `uniformly 
in $r$', in the sense that, for some $M$:
\roster

\item"{I$'$.}" $\sum_{\bx} p^r(\bx) \leq M < \infty  \text{ for all } r$, 

\item"{II$'$.}" for every $\bx \in \Z^{d-1}$, the set $\{t:\P_{\bp^{r}} (0\to 
(\bx, t))>0\}$ is independent of $r$, and has greatest common divisor $1$,
\item"{III$'$.}" there exist an integer $\alpha$ and a 
real $\xi\in(0,1)$ such that
$$
\Sigma_u^r=\sum_{\bx:|\bx|_\oo>u}p^r(\bx)
$$
satisfies $\Sigma_{\a h}^r\le \xi\Sigma_h^r$ for all $h\ge 1$.
\endroster
We impose a further condition, namely:
\roster
\item"{IV$'$.}" for every $\bx\in\Z^{d-1}$, the quantity
$p^r(\bx)$ is a continuous and strictly increasing function of $r$.
\endroster
Suppose that 
$$
\rc =\sup\{r:\P_{\bp^{r}} (0\to \infty) = 0\}
$$
satisfies $0<\rc <1$.
It follows under Conditions I$'$--IV$'$
by the block construction and the standard
 argument quoted in Section 4 that $\theta (\bp^{\rc})=0$. 

\section{Acknowledgements}
PH acknowledges financial support from Trinity College, Cambridge.
GRG thanks the organizers and participants of the Fourth Brazilian
School of Probability for a stimulating meeting in an inspiring location.
The study of infinite-range directed percolation is in part a response
to a question posed to GRG by Gordon Slade.

\Refs
\def\JSP{Journal of Statistical Physics}

\def\CMP{Communications in Mathematical Phy\-sics}

\def\PTRF{Probability Theory and Related Fields}

\def\AP{Annals of Probability}

\def\nref#1\endref{}  % {\ref#1\endref}

\def\MPCPS{Mathematical Proceedings of the Cambridge Philosophical Society}

\ref
\noOf{AKNa}
\by Aizenman, M., Kesten, H., Newman, C.\ M.
\paper Uniqueness of the infinite cluster and related results in
percolation
\inbook Percolation Theory and Ergodic Theory of Infinite Particle
Systems
\ed H.\ Kesten
\publ Springer
\publaddr Berlin--Heidelberg--New York
\yr 1987
\pages 13--20
\bookinfo IMA Volumes in {\it Mathematics and its Applications}
\vol 8
\endref

\ref
\noOf{AN86}
\by Aizenman, M., Newman, C.\ M.
\paper Discontinuity of the percolation density in one-dimensional
$1/|x-y|^2$ percolation models
\jour \CMP
\vol 107
\pages 611--647
\yr 1986
\endref

\ref
\noOf{BGN1}
\by  Barsky, D.\ J.,  Grimmett, G.\ R., Newman, C.\ M.
\paper Dynamic renormalization and continuity of
the percolation transition in orthants
\inbook Spatial Stochastic Processes
\eds K.\ S.\ Alexander and J.\ C.\ Watkins
\yr 1991
\pages 37--55
\publ Birkh\"auser
\publaddr Boston
\endref

\ref
\noOf{BGN2}
\by  Barsky, D.\ J.,  Grimmett, G.\ R., Newman, C.\ M.
\paper Percolation in half spaces: equality of
critical probabilities and continuity of the percolation probability
\jour  Probability Theory and Related Fields
\vol 90 
\yr 1991
\pages 111--148
\endref

\ref
\noOf{BenPP97}
\by Benjamini, I., Pemantle, R., Peres, Y.
\paper Unpredictable paths and percolation
\jour\AP
\vol 26
\pages 1198--1211
\yr 1998
\endref

\ref
\noOf{BeKea}
\by Berg, J.\ van den and Keane, M.
\paper On the continuity of the percolation probability function
\inbook Particle Systems, Random Media and Large Deviations
\ed R.\ T.\ Durrett
\publ AMS
\publaddr Providence, R.\ I.
\bookinfo  Contemporary Mathematics Series
\yr 1984
\vol 26
\pages 61--65
\endref

\ref
\noOf{BG1}
\by Bezuidenhout, C.\ E.,  Grimmett, G.\ R.
\yr 1990
\paper The critical contact process dies out
\jour \AP
\vol 18
\pages 1462--1482
\endref

\ref
\noOf{BK}
\by Burton, R.\ M., Keane, M.
\paper Density and uniqueness in percolation
\yr 1989
\jour \CMP
\vol 121
\pages 501--505
\endref

\ref
\noOf{Cerf0}
\by Cerf, R.
\paper Large deviations for three dimensional supercritical percolation
\yr 2000
\jour Ast\'erisque
\vol 267
\endref

\ref
\noOf{DoyS}
\by Doyle, P.\ G., Snell, E.\ L.
\book Random Walks and Electric Networks
\publ American Mathematical Association
\publaddr Washington, D.\ C
\yr 1984
\bookinfo Carus Mathematical Monograph no.\ 22
\endref

\ref
\noOf{Dur84}
\paper Oriented percolation in two dimensions
\by Durrett, R.\ T.
\jour\AP
\vol 12
\yr 1984
\pages 999--1040
\endref

\ref
\noOf{Dur91}
\by Durrett, R.\ T.
\paper The contact process, 1974--1989
\inbook Mathematics of Random Media, {\rm Blacksburg, VA}
\bookinfo Lectures in Applied Mathematics
\vol 27
\publ AMS
\yr 1991
\pages 1--18
\endref

\ref
\noOf{GanGR}
\by Gandolfi, A., Grimmett, G.\ R., Russo, L.
\paper On the uniqueness of the infinite open cluster in the percolation
model
\jour\CMP
\vol 114
\pages 549--552
\yr 1988
\endref

\ref
\noOf{G99}
\by Grimmett, G.\ R.
\yr 1999
\book Percolation
\bookinfo 2nd edition
\publ Springer 
\publaddr Berlin
\endref

\ref
\noOf{GKM}
\by Grimmett, G.\ R., Keane, M., Marstrand, J.\ M. 
\paper On the connectedness of a random graph
\jour \MPCPS
\vol 96
\pages 151--166
\yr 1984
\endref

\ref
\noOf{GKZ}
\by Grimmett, G.\ R., Kesten, H., Zhang, Y.
\paper Random walk on the infinite cluster of the percolation model
\jour\PTRF
\yr 1993
\vol 96
\pages 33--44
\endref

\ref
\noOf{GrM}
\by Grimmett, G.\ R., Marstrand, J.\ M. 
\paper The supercritical phase of percolation is well
behaved
\jour Proceedings of the Royal Society (London), Series A
\vol  430
\yr 1990
\pages 439--457
\endref 

\ref\noOf{HM98}
\by H\"aggstr\"om, O.,  Mossel, E.
\paper Nearest-neighbor walks with low predictability profile and
percolation in $2+\epsilon$ dimensions
\jour\AP
\vol  26 
\yr 1998
\pages 1212--1231
\endref

\ref\noOf{HP99}
\by H\"aggstr\"om, O., Peres, Y.
\paper Monotonicity of uniqueness for percolation
on Cayley graphs: all infinite clusters are born simultaneously
\jour\PTRF
\vol  113 
\yr 1999
\pages 273--285
\endref

\ref\noOf{HoM}
\by Hoffman, C.\ and Mossel, E.
\paper Energy of flows on percolation clusters
\jour Potential Analysis
\yr 1999
\toappear
\endref

\ref
\noOf{vdHS00}
\by Hofstad, R.\ van der, Slade, G.
\paper A generalised inductive approach to the lace expansion
\toappear
\yr 2000
\endref

\ref
\noOf{Kes81}
\by Kesten, H.
\paper Analyticity properties and power law estimates in percolation
theory
\jour\JSP
\vol 25
\pages 717--756
\yr 1981
\endref

\ref
\noOf{Lig}
\by Liggett, T.\ M.
\book Interacting Particle Systems
\publ Springer
\publaddr Berlin
\yr 1985
\endref

\ref
\noOf{Lig99}
\by Liggett, T.\ M.
\book Stochastic Interacting Systems: Contact, Voter
and Exclusion Processes
\publ Springer
\publaddr Berlin
\yr 1999
\endref

\ref
\noOf{LSS}
\by Liggett, T.\ M., Schonmann, R.\ H., Stacey, A.
\paper Domination by product measures 
\yr 1997
\jour\AP
\vol 25
\pages 71--95
\endref

\ref
\noOf{MeeS}
\by Meester, R., Steif, J.
\paper On the continuity of the critical value for long range percolation 
in the exponential case
\jour\CMP
\vol 180
\pages 483--504
\yr 1996
\endref 

\ref
\noOf{Pis95}
\by Pisztora, A.
\paper Surface order large deviations for Ising, Potts and percolation models
\jour \PTRF
\vol 104
\pages 427--466
\yr 1996
\endref

\ref
\noOf{Sch99}
\by Schonmann, R.\ H.
\paper Stability of infinite clusters in supercritical percolation
\jour \PTRF
\vol 113
\yr 1999
\pages 287--300
\endref

\ref
\noOf{Sl00}
\by Slade, G.
\paper Lattice trees, percolation and super-Brownian motion
\toappear
\yr 2000
\endref

\ref
\noOf{Soa}
\by Soardi, P.\ M.
\book Potential Theory on Infinite Networks
\publ Springer
\publaddr Berlin
\yr 1994
\endref

\endRefs
\enddocument